\documentclass[12pt, dvips, twoside]{amsart}
\usepackage{amssymb,latexsym,bbm,graphicx,epsfig,epic,eepic,oldgerm}
\usepackage{psfrag}

\theoremstyle{plain}

\newcommand{\proofend}{\hspace*{\fill} $\Box$\\}

\def\s{\smallskip}
\def\m{\medskip}
\def\eps{\varepsilon}

\def\HZ{\operatorname{HZ}}

\def\eps{\epsilon}

\def\gf{\varphi}

\def\go{\omega}

\def\ch{{\mathcal H}}

\def\cl{{\mathcal L}}

\def\cp{{\mathcal P}}

\def\NN{\mathbbm{N}}

\def\RR{\mathbbm{R}}

\def\pp{\partial}

\def\ra{\rightarrow}

\def\ni{\noindent}
\def\b{\bigskip}
\def\m{\medskip}

\def\proof{\noindent {\it Proof. \;}}

\begin{document}

\begin{titlepage}
\title[A refinement of the Hofer--Zehnder theorem]{A refinement of the Hofer--
Zehnder theorem on the existence of closed trajectories near a hypersurface}
$   $ \\
$   $ \\

\author{Leonardo Macarini$^1$}
\address{(L.\ Macarini) Instituto de Matem\'atica Pura e Aplicada - IMPA\\
         Estrada Dona Castorina, 110 - Jardim Bot\^anico\\
         22460-320 Rio de Janeiro RJ\\
         Brazil}
\email{leonardo@impa.br}
\author{Felix Schlenk$^2$}
\address{(F.\ Schlenk) ETH Z\"urich, CH-8092 Z\"urich, Switzerland}
\email{schlenk@math.ethz.ch}

\thanks{$^1$ Partially supported by CNPq-PROFIX, Brazil}
\thanks{$^2$ Supported by the von Roll Research Foundation}

\date{\today}

\end{titlepage}

\begin{abstract}
The Hofer--Zehnder theorem states that almost every hypersurface in a
thickening of a hypersurface $S$ in a symplectic manifold $(M,\omega)$
carries a closed characteristic provided that $S$ bounds a compact
submanifold and $(M,\omega)$ has finite capacity.
We show that it is enough to assume that the thickening of $S$ has
finite capacity.
\end{abstract}

\maketitle

\ni
We consider a smooth symplectic manifold $(M,\go)$.
A {\it hypersurface} $S$ in $M$
is a smooth compact connected orientable codimension
$1$ submanifold of $M \setminus \pp M$ without boundary.
A characteristic on $S$ is an embedded circle in $S$ all
of whose tangent lines belong to the
distinguished line bundle 
\[
\cl_S \,=\, \left\{ (x, \xi) \in TS \mid \go(\xi, \eta) =0 
\text{ for all } \eta \in T_x S \right\} .
\]
We denote by $\cp (S)$ the set of closed characteristics on $S$.
Examples show that $\cp (S)$ can be empty, see \cite{Gi, GG}.
We therefore follow \cite{HZ} and consider parametrized neighbourhoods of $S$.
Since $S$ is orientable, there exists 
an open neighbourhood $I$ of $0$ and a smooth diffeomorphism 
\[
\psi \colon S \times I \,\ra\, U \subset M
\]
such that $\psi (x,0) =x$ for $x \in S$.
We call $\psi$ a {\it thickening of $S$}, and we
abbreviate $S_\eps = \psi \left( S \times \left\{ \eps \right\} \right)$.

Given an open subset $U \subset M$ we consider the function space
$\ch(U)$ of smooth functions $H \colon U \ra [0,\max H]$ such that
\begin{itemize}
\item[$\bullet$]
$H |_V = 0$ for some nonempty open set $V \subset U$;
\item[$\bullet$]
$H |_{U\setminus K} = \max H$ for some compact set $K \subset U$.
\end{itemize}
We say that $H \in \ch(U)$ is {\it $\HZ$-admissible}\, if the flow
$\gf_H^t$ has no non-constant $T$-periodic orbit with period $T \le 1$,
and we set
\[
\ch_{\HZ} (U) \,=\, \left\{ H \in \ch(U) \mid H \text{ is
HZ-admissible} \right\} .
\]
The Hofer--Zehnder capacity of $U$ is defined as
\[
c_{\HZ} (U) \,=\, 
      \sup \left\{ \max H \mid H \in \ch_{\HZ} (U) \right\} .
\]

It has been shown in \cite[Sections~4.1 and 4.2]{HZ} that for any thickening 
$\psi \colon S \times I \ra U \subset M$ for which 
$c_{\HZ} \left( U \right) < \infty$,
the set $\left\{ \eps \in I \mid \cp \left( S_{\eps} \right) \neq \emptyset \right\}$ is dense in $I$, 
and that
\[
\mu \left\{ \eps \in I \mid \cp \left( S_{\eps} \right) 
\neq \emptyset \right\} \,=\, \mu (I)
\]
if $S$ bounds a compact submanifold of $M$ and $c_{\HZ} (M) < \infty$.
Here, $\mu$ denotes Lebesgue measure. 
In this note we prove

\b \ni
{\bf Theorem~1.}
{\it
For any thickening $\psi \colon S \times I \ra U \subset M$  we have
\[
\mu \left\{ \eps \in I \mid \cp \left( S_\eps \right)  
\neq \emptyset \right\} \,=\, \mu (I) 
\]
provided that $c_{\HZ} (U) < \infty$.
}

\b
\proof
Consider a thickening $\psi \colon S \times I \ra U \subset M$.
We can assume that $I = \;]-1,1[$.
Let $G = \left\{ \eps \in \;]-1,1[ \; \mid \cp \left( S_{\eps} \right) 
\neq \emptyset \right\}$ be the parameters of the good hypersurfaces.

\m
\ni
{\bf Step 1.} {\it The set $G$ is measurable.}

\s \ni
\proof
We define the smooth function $K \colon U \ra \;]-1,1[$ by
\[
K(z) = \eps \quad \text{if }\, z \in S_{\eps} .
\]
The set $\cp_K (\eps)$ of closed orbits of the Hamiltonian flow $\gf_K^t$ 
of $K$ on $S_\eps$ corresponds to $\cp \left( S_\eps \right)$.
For each closed orbit $x$ of $\gf_K^t$ we denote by $T(x)$ 
its period defined as
\[
T(x) \,=\, \min \left\{ t>0 \mid \gf_K^t (p) = p \right\} 
\]
where $p$ is any point on $x$.
For each $n \in \NN$ we set 
\[
G_n \,=\, \left\{ \eps \in\; ]-1,1[\; \mid \text{ there exists } 
              x \in \cp_K (\eps) \text{ such that } T(x) \le n \right\} .
\]
It follows from the Arzel\`a--Ascoli Theorem that
$G_n$ is a closed subset of $]-1,1[$, see 
\cite[page 109, Proposition 1]{HZ}
The set $G = \bigcup_{n\ge 1} G_n$ is thus measurable.

\m
\ni
{\bf Step 2.} {\it The set $G$ has full measure.}

\s \ni
\proof
Applying the argument in the proof of Theorem~3 in Section~4.2
of \cite{HZ} to Hamiltonian functions $H \colon U \ra \RR$ 
as in the figure below
such that
\[
H(x) \,=\, f \left( \left| K(x) \right| \right) \quad \text{if }\, x \in
S_{-t} \cup S_t,\; \eps \le t \le 1 ,
\]
we see that $\cp \left( S_{-\eps} \right) \neq \emptyset$ or
$\cp \left( S_\eps \right) \neq \emptyset$ whenever the function 
$C(t) \colon [0,1[\; \ra \RR$, 
$t \mapsto c_{\HZ} \left( \psi \left( S \times ]-t,t[ \right) \right)$ 
is Lipschitz continuous at $\eps$.
Since $C(t)$ is monotone increasing, 
it is differentiable almost everywhere, 
and thus Lipschitz continuous almost everywhere.
Therefore,
\begin{equation}  \label{e:m1}
\mu \left\{ 
    \eps \in [0,1[\; \mid \cp \left( S_{-\eps} \right) \neq \emptyset
    \text{ or }         \cp \left( S_{\eps} \right) \neq \emptyset 
    \right\} \,=\, 1 .
\end{equation}
\begin{figure}[h] 
 \begin{center}
  \psfrag{S}{$S$}
  \psfrag{e}{$S_\eps$}
  \psfrag{-e}{$S_{-\eps}$}
  \psfrag{es}{$S_{\eps^*}$}
  \psfrag{-es}{$S_{-\eps^*}$}
  \psfrag{Ce}{$C(\eps)$}
  \psfrag{Ces}{$C(\eps^*)$}
  \leavevmode\epsfbox{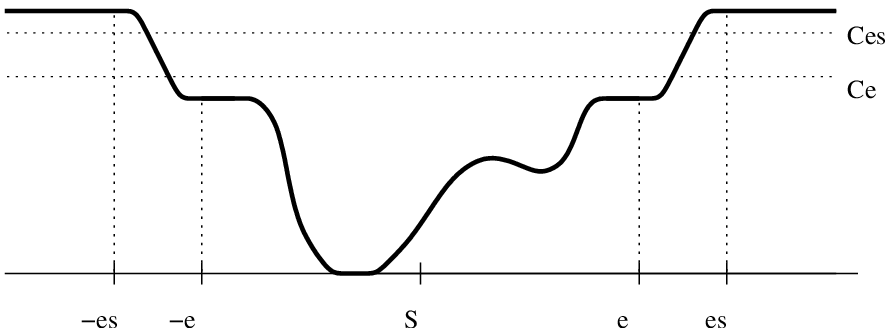}
 \end{center}
 \label{hoz}
\end{figure}
%
%

\ni
Let $B = \;]-1,1[ \,\setminus G$ be the parameters of the bad hypersurfaces.
Arguing by contradiction we assume that $\mu (B) >0$.
Choose an interval $]a,b[\; \subset \;]-1,1[$ such that 
\begin{equation}  \label{e:23}
\mu \left( ]a,b[\, \cap \:\!B \right) \,\ge\, \tfrac 23 (b-a) .
\end{equation} 
Applying \eqref{e:m1} to $]a,b[$ instead of $]-1,1[$ we see that
\begin{equation}  \label{e:12}
\mu \left( ]a,b[\, \cap \:\!G \right) \,\ge\, \tfrac 12 (b-a) .
\end{equation} 
Combining \eqref{e:23} with \eqref{e:12} we find that
\begin{eqnarray*}
b-a &=& \mu (]a,b[) \\
    &=& \mu \left( ]a,b[\, \cap \:\!B \right) + \mu \left( ]a,b[\, \cap \:\!G \right) \\
    &\ge& \tfrac 23 (b-a) + \tfrac 12 (b-a) \\
    &>& b-a .
\end{eqnarray*}
This contradiction shows that $\mu (B) =0$, and so 
$\mu (G) = \mu \left( ]-1,1[ \right) $.
\proofend

It is often of interest to find closed characteristics in a restricted set
of homotopy classes. 
We thus fix a subset $\Gamma$ of the fundamental group $\pi_1(S)$,
and given a thickening $\psi \colon S \times I \ra U \subset M$,
we denote by $\cp^\Gamma \left( S_\eps \right)$ the set of 
closed characteristics on $S_\eps$ representing an element of 
$\Gamma \subset \pi_1 \left( S_\eps \right) = \pi_1(S)$.
We can assume that $I = \;]-1,1[$,
and for each $\eps \in \;]0,1[$ we set 
$U_\eps = \psi \left( S \times ]-\eps, \eps [\right)$.
We say that $H \in \ch \left( U_\eps \right)$ is 
{\it $\HZ^\Gamma$-admissible}\, 
if the flow $\gf_H^t$ has no non-constant $T$-periodic orbit 
with period $T \le 1$ which represents an element of 
$\Gamma \subset \pi_1 \left( U_\eps \right) = \pi_1(S)$,
and we set
\[
\ch_{\HZ}^\Gamma \left( U_\eps \right) \,=\, 
 \left\{ H \in \ch \left( U_\eps \right) 
   \mid H \text{ is $\HZ^\Gamma$-admissible} \right\}.
\]
The $\Gamma$-sensitive Hofer--Zehnder capacity of $U_\eps$ 
is defined as
\[
c_{\HZ}^\Gamma \left( U_\eps \right) \,=\, 
\sup \left\{ \max H \mid H \in
\ch_{\HZ}^\Gamma \left( U_\eps \right) \right\} .
\]
Notice that the map $\eps \mapsto c_{\HZ}^\Gamma \left( U_\eps \right)$ 
is monotone increasing.
Repeating the above proof with 
$c_{\HZ}$ replaced by $c_{\HZ}^\Gamma$ we find
the following refinement of Theorem~1.

\b \ni
{\bf Theorem~2.}
{\it
For any thickening $\psi \colon S \times I \ra U \subset M$ of a
hypersurface $S$ in $(M,\go)$ we have
\[
\mu \left\{ \eps \in I \mid \cp^\Gamma \left( S_\eps \right)  
\neq \emptyset \right\} \,=\, \mu (I) 
\]
provided that $c_{\HZ}^\Gamma (U) < \infty$.
}

\b
\ni
While a version of Theorem~1 is applied in \cite{FS},
Theorem~2 is used in \cite{Ma} to show that the flow 
describing the dynamics of a charge in a rational symplectic magnetic 
field has a closed orbit on almost every small energy level.

\b
\ni
{\bf Acknowledgements.}
This note was written during the second authors stay at Stony Brook in
April 2003. 
He wishes to thank Ely Kerman and Dusa Mc\;\!Duff 
for their warm hospitality.


\begin{thebibliography}{99}

\bibitem{FS}
U.\ Frauenfelder and F.\ Schlenk.
Hamiltonian dynamics on convex symplectic manifolds.
math.SG/0303282.
    

\bibitem{Gi}
V.\ Ginzburg.
Hamiltonian dynamical systems without periodic orbits. 
Northern California Symplectic Geometry Seminar, 35--48, 
{\it Amer. Math. Soc. Transl. Ser. 2}, {\bf 196}, 
Amer. Math. Soc., Providence, RI, 1999. 


\bibitem{GG}
V.\ Ginzburg and B.\ G\"urel.
Relative Hofer--Zehnder capacity and periodic orbits in twisted 
cotangent bundles.
math.DG/0301073. 


\bibitem{HZ} 
H.\ Hofer and E.\ Zehnder. 
{\it Symplectic Invariants and Hamiltonian Dynamics}. 
Birkh\"auser, Basel, 1994.


\bibitem{Ma}
L.\ Macarini.
Hofer--Zehnder semicapacity of cotangent bundles and symplectic submanifolds.
math.SG/0303230.



\end{thebibliography}
\enddocument